\documentclass{amsart}
\usepackage{amssymb}
\usepackage{graphicx}
\usepackage{url}
 \newtheorem{theorem}{Theorem}
 \newtheorem{lemma}[theorem]{Lemma}
 \newtheorem{proposition}[theorem]{Proposition}
 \newtheorem{corollary}[theorem]{Corollary}
 \newtheorem{remark}[theorem]{Remark}
 \newtheorem{example}[theorem]{Example}
 \newtheorem{definition}[theorem]{Definition}
 \newtheorem{conjecture}[theorem]{Conjecture}

 \newtheorem{question}[theorem]{Question}

\newcommand{\bpr}{\begin{proof}}
\newcommand{\epr}{\end{proof}}
\newcommand{\beq}{\begin{equation}}
\newcommand{\eeq}{\end{equation}}
\newcommand{\bThm}{\begin{theorem}}
\newcommand{\eThm}{\end{theorem}}
\newcommand{\blem}{\begin{lemma}}
\newcommand{\elem}{\end{lemma}}
\newcommand{\bpro}{\begin{proposition}}
\newcommand{\epro}{\end{proposition}}
\newcommand{\bcor}{\begin{corollary}}
\newcommand{\ecor}{\end{corollary}}
\newcommand{\brem}{\begin{remark}}
\newcommand{\erem}{\end{remark}}
\newcommand{\bexa}{\begin{example}}
\newcommand{\eexa}{\end{example}}
\newcommand{\bdf}{\begin{definition}}
\newcommand{\edf}{\end{definition}}
\newcommand{\bcon}{\begin{conjecture}}
\newcommand{\econ}{\end{conjecture}}
\newcommand{\bque}{\begin{question}}
\newcommand{\eque}{\end{question}}

\newcommand{\p}{\partial}

\newcommand{\comment}[1]{}

\title{Ribbon $n$-knots with isomorphic quandles}
\author{Blake K. Winter}
\thanks{The author wishes to acknowledge Adam Sikora, William Menasco, Jason Manning, and Sam Nelson, for their recommendations on the writing of this paper.}
\address{}

\date{}

\begin{document}
\thispagestyle{empty}

\begin{abstract}
Let $K, K'$ be ribbon knottings of $n$-spheres with $1$-handles in $S^{n+2}$, $n\geq 2$. We show that if the knot quandles of these knots are isomorphic, then the ribbon knottings are stably equivalent, in the sense of Nakanishi and Nakagawa, after taking a finite number of connected sums with trivially embedded copies of $S^{n-1}\times S^{1}$.
\end{abstract}

\maketitle

%
\section{Introduction}
%

The \emph{fundamental quandle} of a knot or link is a powerful invariant for classical links in $S^3$, \cite{DJ, Mat}. It is also an invariant for links in any dimension. In this paper, we will consider ribbon $n$-links, $n \geq 2$, in $S^{n+2}$, under an equivalence relation wherein two links are considered to be equivalent if they are isotopic after taking the connected sum with a finite number of trivially embedded copies of $S^1 \times S^{n-1}$. We will show that the quandle is a complete invariant for this equivalence relation: two ribbon links are equivalent iff they have isomorphic quandles. Furthermore, we will show that the isotopy takes on the special form of a \emph{stable equivalence}, as defined in \cite{Naka, NakaNaka}.

\section{Ribbon Links}


Here we will review the definition of ribbon links. See \cite{CKS, Naka, NakaNaka} for further discussion of ribbon links in arbitrary dimensions.

\begin{remark}
Throughout, we will use the notation $D^n$ to represent a closed $n$-disk.
\end{remark}

Intuitively, a ribbon presentation is a collection of embedded $(n+1)$-disks, joined together by $1$-handles that are attached to the boundaries of the disks, and which are allowed to pass through the interiors of the disks transversely. We will make this precise: 

Let $M$ be an $(n+2)$-manifold. Let $B$ be a set of disjoint oriented $(n+1)$-disks $B=\{B_{i}:D^{n+1}\rightarrow M\}$. Let $H$ be a set of oriented $1$-handles $H_{j}: D^1\times D^n\rightarrow M$, with the following conditions:
\begin{enumerate}
	\item $H_j(\{ 0,1\}\times D^n)$ must be embedded in $\cup_{B_i\in B}B_i(\partial D^{n+1})$, with compatible orientations.
	\item For all $i, j$, the components of $H_j((0,1)\times D^n)\cap B_i(D^{n+1})$ are a finite embedded collection of $n$-disks embedded in $H_j((0,1)\times D^n)$, such that each component separates $H_j(\{0\}\times D^n)$ from $H_j(\{1\}\times D^n)$.
\end{enumerate}

This implies that $H_j((0,1)\times D^n)\cap B_i(D^{n+1})$ is contained in the interior of the image of $B_i$. Such intersections between a handle and a base are called \emph{ribbon intersections} or \emph{ribbon singularities}. A handle can intersect any number of bases in ribbon intersections, and it can intersect the same base multiple times (or zero). We will in general identify the maps $B_i$ and $H_j$ with their images.
\begin{definition}
The disks in $B$ are called \emph{bases}, while the elements of $H$ are called \emph{(fusion) handles} or \emph{bands}.
\end{definition}
\begin{definition}
Given a set of bases and handles, the union of all the bases and handles together forms an $(n+1)$ immersed manifold, which is called a \emph{ribbon solid}.
\end{definition}
Any immersed $(n+1)$-manifold which can be expressed in this form will also be called a ribbon solid. Note that a given ribbon solid might be able to be broken up into bases and handles in many different ways.

\begin{definition}
Let $(B, H)$ be collections of bases and handles. Let $K$ be the link defined as $K=({\cup}_{B_i\in B}(\partial B_{i}-{\cup}_{H_j\in H}H_{j}))\cup({\cup}_{H_j\in H}H_{j}(D^1\times\partial D^n))$. Then $K$ is called a \emph{ribbon link}, and the pair $(B, H)$ are a \emph{ribbon presentation} for $K$.
\end{definition}

In other words, $K$ is the closure of the boundary of the ribbon solid defined by the union of the bases and handles in the pair $(B, H)$. We identify two ribbon presentations, $(B, H)$ and $(B', H')$, of $K$ if there is a path between them in the space of ribbon presentations of $K.$ This is equivalent to requiring that $|B|=|B'|$, $|H|=|H'|$, and that there exists an ambient isotopy $f_t$ of $M$, $t\in [0,1]$, such that $f_0$ is the identity, and for all $B_i\in B, H_j\in H, B'_i\in B', H'_j\in H'$, $f_1\circ B_i=B'_i$, and $f_1\circ H_j=H'_j$.

Observe that since we require that the handles be glued to the bases in an orientation-preserving manner, $K$ will be orientable. Note also that for ribbon $1$-links, the components of a ribbon link may not be ribbon knots themselves. However, for $n$-links with $n\geq 2$, each component of a ribbon link will automatically be a ribbon knot.

Given a ribbon link $K$, there are in general multiple ribbon presentations for $K$. There are some distinct notions of equivalence for ribbon presentations. We will, for the moment, restrict our attention to knots; the generalization to links is straightforward. 

\begin{definition}
Let $K$ be a ribbon knot with ribbon presentations $(B, H)$ and $(B, H')$. These presentations are called \emph{simply equivalent} if $|H|=|H'|$, $|B|=|B'|$, and there exists an isotopy $f_t$, $t\in [0,1]$, of $M$, such that $f_0$ is the identity, and for all $B_i\in B, H_j\in H, B'_i\in B', H'_j\in H'$, $f_1\circ B_i|_{\partial D^{n+1}}=B'_i|_{\partial D^{n+1}}$, and $f_1\circ H_j=H'_j$.

\end{definition}

\begin{figure}
		\centering
			\includegraphics[scale=0.3]{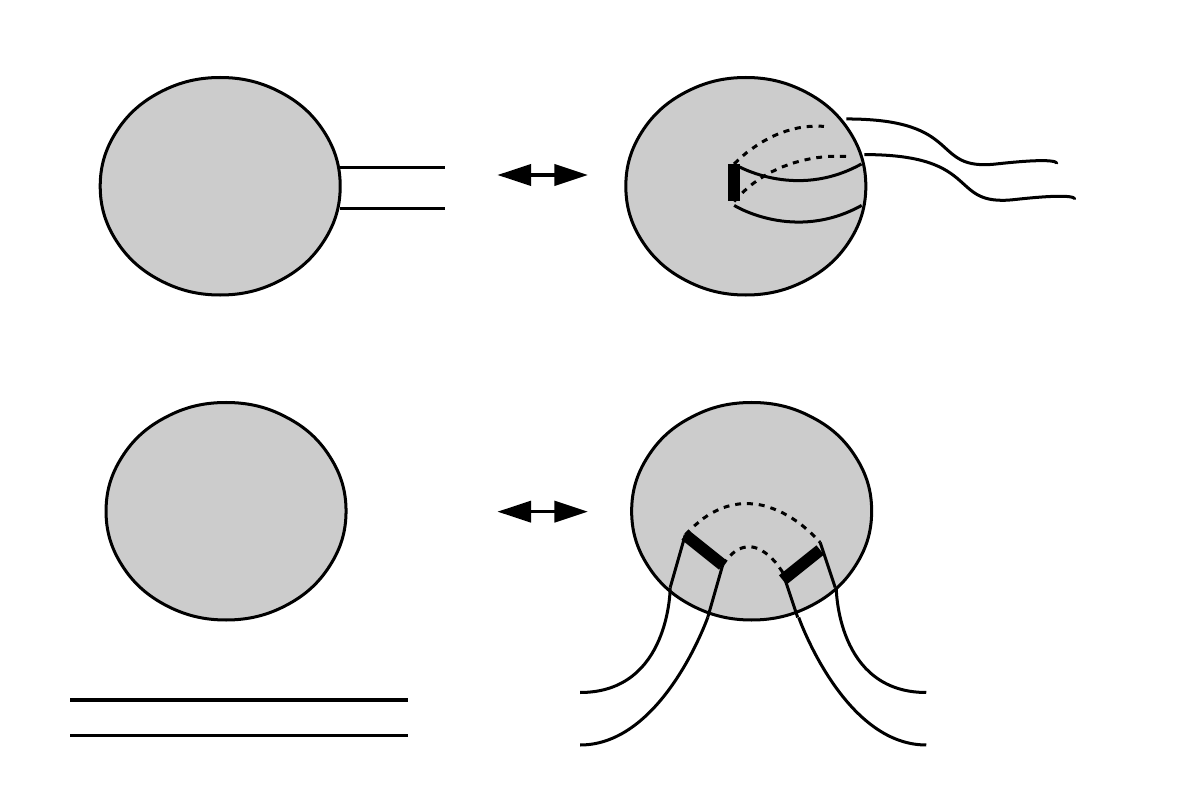}
		\caption{Two simply equivalent ribbon presentations with distinct ribbon intersections along a single handle. Here, the bases are indicated with shading; ribbon intersections are shown by thickened lines.}
		\label{ribbonstable2}
\end{figure}

That is, the isotopy maps handles to handles, and maps the boundaries of bases to the boundaries of bases. It need not map the interiors of the bases to one another, however. Thus, simple equivalence may create new intersections of handles with the interiors of bases, as illustrated in Fig. \ref{ribbonstable2}.

\begin{lemma}
For $n\geq 2$ a ribbon presentation $(B,H)$ in a simply-connected manifold $M$ is determined uniquely up to simple equivalence by specifying the following ribbon data:\\
(a) the cardinality of $B$ and,\\
(b) for each $H_{j}\in H$, a specification of the base to which $H_j (\{ 0 \} \times D^n)$ is attached, a specification of the base to which $H_j (\{ 1 \} \times D^n)$ is attached, and a sequence of bases, $B_{j1},...,B_{jk},$
through which the handle $H_j$ passes through, up to reversal of that sequence.
\end{lemma}

\bpr Clearly there is a ribbon presentation fulfilling every possible ribbon data as above.
We need to show that such ribbon presentation is unique, up to simple equivalence. Note that for for each $j$, the sequence $B_{j1},...,B_{jk},$ determines the homotopy class of the handle $H_j$ in $M-\cup_i \p B_i$ uniquely. Since $n\geq 2$ and the core of the handle is $1$-dimensional, the homotopy class of the handle determines its isotopy class as well.
\epr

In $n=1$ this is of course not the case, since the homotopy class of the handle does not determine its isotopy class in $M-\cup_i \p B_i$. In addition, one must specify a cyclic ordering  of handles glued to a base for $n=1$.

There is a notion of equivalence between ribbon presentations called \emph{stable equivalence}, defined by Nakanishi, \cite{Naka}. 

\begin{figure}
		\centering
			\includegraphics[scale=0.6]{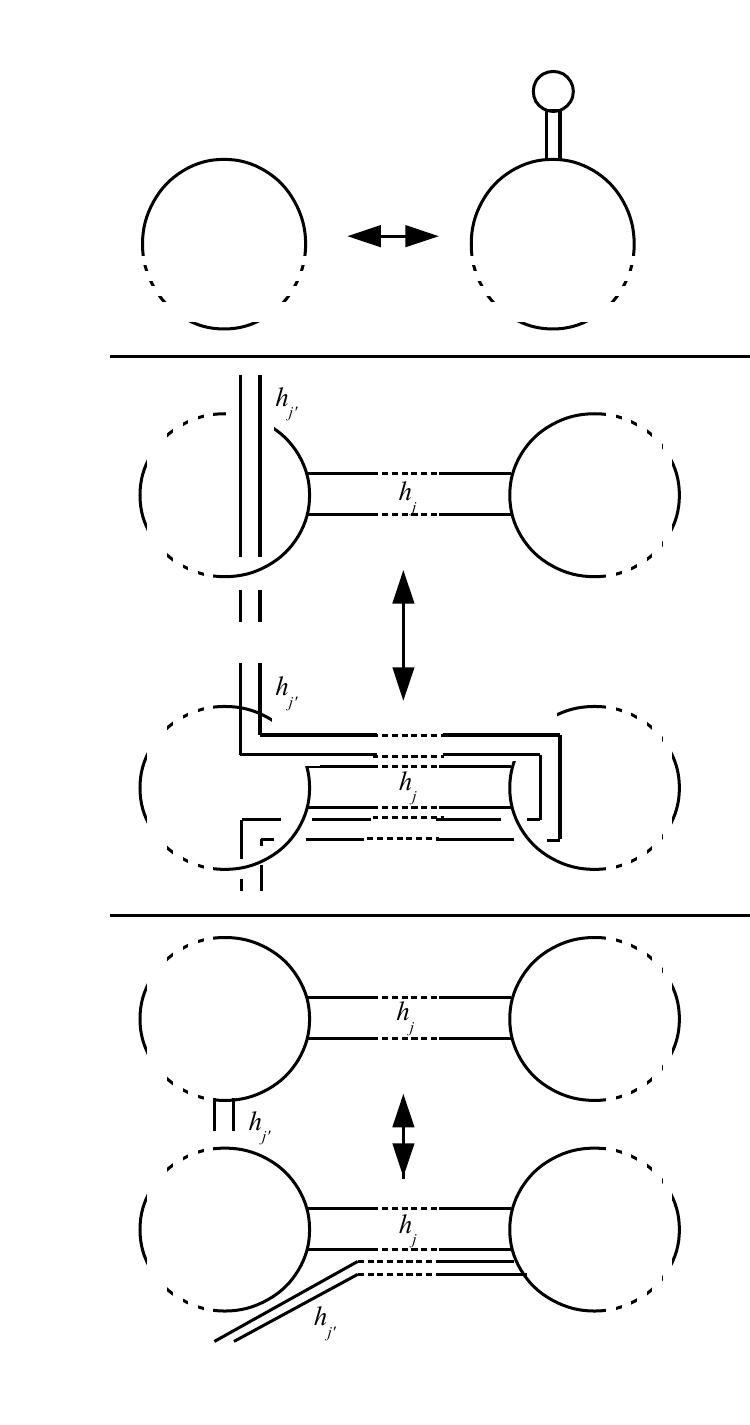}
		\caption{The three moves permitted for ribbon presentations under stable equivalence.}
		\label{ribbonstable}
\end{figure}

\begin{definition}
Two ribbon presentations are \emph{stably equivalent} if they are related by simple equivalence, together with the following three changes, which are illustrated in Fig. \ref{ribbonstable}:

\begin{enumerate}
	\item Add a new base $B_{0}$ with a handle $H_{0}$ starting on $B_{0}$ and terminating on any other base, without passing through any bases.
	\item \emph{Handle pass}: Isotope $H_{j}$ through $H_{j'}([0,1]\times (D^n-\partial D^n))$, leaving $H_j|_{\{ 0,1\}\times D^n}$ fixed.
	\item \emph{Handle slide}: If $H_{j}$ connects $B_{i}$ with $B_{i'}$, and $H_{j'}$ has an end on $B_{i}$, we may slide that end along $H_{j}$ over to $B_{i'}$.
\end{enumerate}
\end{definition}

Observe that for a knotted $n$-sphere, $|H|=|B|-1$ by Euler characteristic arguments. It is easy to check that for $n\geq 2$, a ribbon $n$-knot $K$ will have the homeomorphism type of an $n$-sphere with $|H|-|B|+1$ $1$-handles attached.

Given a ribbon $n$-knot $K$, it is an open question whether any two ribbon presentations for $K$ are related by stable equivalence for $n\geq 2$. For the case $n=1$, Nakanishi and Nakagawa provide a family of knots, each of which admits multiple ribbon presentations which are not stably equivalent, \cite{NakaNaka, Naka}.

\section{Quandles}\label{QaB}

\begin{definition}
A \emph{quandle}, \cite{DJ}, (also called a \emph{distributive groupoid} in \cite{Mat}) $(Q, *)$ consists of a set $Q$ and a binary operation $*$ on $Q$ such that for all $a, b, c \in Q$, the following equalities hold.

\begin{enumerate}
	\item $a*a=a$.

\item There is a unique $x\in Q$ such that $x*b=a$.

\item $(a*b)*c=(a*c)*(b*c)$.

\end{enumerate}
\end{definition}

When the operation $*$ is implied, we may denote the quandle by $Q$, suppressing the operation.

In general, it is common to denote the $*$ operation in any quandle $Q$ by using the conjugation notation. For any quandle we may define $a^{b}=a*b$. Analogously, we will denote the element $x$ stipulated by property (2) by $a^{\bar{b}}$, i.e. $(a^b)^{\bar{b}}=a$. By convention, we interpret $a^{bc}=(a^b)^c$.

Quandles as knot invariants were introduced in \cite{DJ, Mat}.
Here we recall the definition of a quandle associated to a knot, as given by Joyce for classical knots in \cite{DJ}. The geometric definition holds for any codimension-$2$ knot, however.


Let $L$ be a link in $M$, with $L$ and $M$ assumed to be orientable, and let $N(L)$ be an open tubular neighborhood of $L$. Then $M-N(L)$ is a manifold with boundary containing $\partial \overline{N(L)}$. The fundamental group of this space is called the \emph{link group}. We will follow the convention that the product of two homotopy classes of paths goes from right to left, i.e. $\gamma' \gamma$ is the homotopy class of a path which follows $\gamma$ and then $\gamma'$. This will make it somewhat easier to discuss the action of the link group on the link quandle.

A meridian $m$ of $L$ is the boundary of the disk fibre of the tubular neighborhood of $L$ considered as a disk bundle over $L$. Note that the orientations of $M$ and of $L$ determine an orientation of $m$. A meridian is defined uniquely up to conjugation.

\begin{figure}[htbp]
\begin{centering}
\includegraphics[scale=.8]{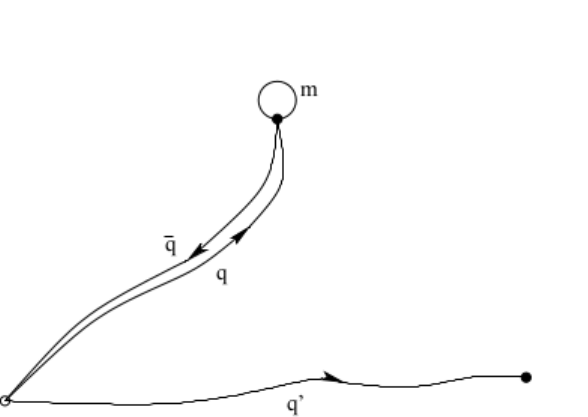}
\caption{The quandle operation for the fundamental quandle of a link $L$ in any space $M$. The
open circle indicates the basepoint of $M$, while the black circles represent points on the boundary of a tubular neighborhood of $L$.}
\label{qop}
\end{centering}
\end{figure}

We define the quandle $Q(M,L)$ associated to the pair $(M,L)$ and a chosen base point $x\in M-N(L)$ as follows. The elements of $Q(M,L)$ are homotopy classes of paths in $M-N(L)$ that start at $x$ and end in $\partial \overline{N(L)}$, where homotopies are required to preserve these two conditions. The quandle operation, $[q]^{[q']}$, between two elements $[q], [q']\in Q(M,L)$ is defined as the homotopy class of $q\overline{q'}mq'$, where the overbar indicates that the path is followed in reverse and $m$ is the meridian of $L$ passing through the endpoint of $q$, cf.  \ref{qop}. It is a straightforward exercise to show that this binary operation is well-defined and it satisfies the definition of a quandle operation. We call $Q(M,L)$ the \emph{link quandle} of $L$ (or the knot quandle if $L$ is a knot).

Note that there is a right action of $\pi_{1}(M-N(L))$ on $Q(M,L)$: $[q]g$ is the equivalence class including the homotopy class of the path $qg$.

The following theorem was proved for classical knots in \cite{DJ}. The arguments in \cite{DJ} are easily generalized to any knots, as is done in \cite{BKW}. As we will not be making use of this result, we will skip the proof here, referring interested readers to those two sources.

\begin{theorem}\label{main}
For knots, $Q(M,K)$ is determined by $\pi_{1}(M-K)$ with its peripheral structure.
\end{theorem}

%
\section{Ribbon $n$-Knots with Isomorphic Quandles}
\label{s-Ribbon-Isom-Q}
%

Let us say that in a ribbon presentation $(B, H)$ a handle $h\in H$ is a \emph{trivial handle} if both ends are attached to the same base and $h$ can be isotoped so it does not meet any bases except at its ends. In other words, $h$ does not pass through any bases.

\begin{definition}
Let $K, K'$ be $n$-knots in $S^{n+2}$, and let $T=S^1\times S^{n-1}$. Let $T_n$ indicate the connected sum of $n$ copies of $T$ that are trivially embedded, that is, which bound embedded copies of $S^1\times D^{n}$. We will say that $K$ and $K'$ are \emph{weakly equivalent} if $K\# T_n$ is isotopic to $K'\# T_m$.
\end{definition}

Note that in \cite{Liv}, weak equivalence is called ``stable equivalence,'' as gluing on copies of $T$ may be thought of as a stabilization of the knot. Because this term is also used to refer stable equivalence of ribbon presentations, we are adopting the term ``weak equivalence'' instead to avoid any confusion.

\begin{definition}
Let $(B, H), (B', H')$ be ribbon presentations for $n$-links, $n\geq 2$. We will say that they are \emph{$k$-weakly stably equivalent} if they are stably equivalent after attaching $k$ trivial handles, and that they are \emph{weakly stably equivalent} if they are $k$-weakly stably equivalent for some finite $k$.
\end{definition}

Given a ribbon presentation $(B, H)$, the bases $b_i$ in $B$ are a collection of embedded $(n+1)$ disks. We will say that a path in $S^{n+2}-K$ \emph{passes through} a base $b_i\in B$ if it intersects the interior of $b_i$ transversely. Note that in general position, a path can pass through a base positively or negatively, depending on whether it passes through it in the same direction as the normal vector or in the opposite direction.

\begin{remark}
To specify a handle $h\in H$, it suffices to specify the starting base to which $\{ 0\} \times D^n$ is glued, the signed intersections of the handle in the order it passes through bases, and the terminal base to which $\{ 1\} \times D^n$ is glued. Any two handles meeting these conditions are smoothly isotopic and thus equivalent.
\end{remark}

\begin{theorem}
Let $K$ be a ribbon knot in $S^{n+2}$, $n\geq 2$, of arbitrary genus. Suppose that $K$ admits a base and handle presentation $(B, H)$ with $k$ bases, and $\pi_{1}(S^{n+2}-K)\cong \mathbb{Z}$. Then $(B, H)$ is $(k-1)$-weakly stably equivalent to an unknotted ribbon knot.
\end{theorem}
\bpr Let $h_{2},..., h_{k}$ be $(k-1)$ new trivial handles glued to $K$, with all their ends connected to base $B_{1}$. Perform handle slides so that $h_{i}$ starts on $B_{1}$ and terminates on $B_{i}$, twisting them at the end so that they pass through bases equally many times positively as negatively. Let $\gamma_{i}$ be a curve which travels to $B_{1}$ without passing through any bases, then follows the handle $h_{i}$ to $B_{i}$ (traveling through the same bases in the same order) and returns to the basepoint for the fundamental group without passing through any additional bases. Since $\pi_{1}(S^{n+2}-K)\cong \mathbb{Z}$, and since by construction $\gamma$ is homologically trivial, $\gamma$ bounds a disk in $S^{n+2}-(K\cup h_{2}\cup ... \cup h_{k})$. Note in particular that this disk does not pass through the handles $h_2, ..., h_k$. This disk may be taken to be smooth and immersed by various approximation theorems, \cite[Theorem 10.16, Theorem 10.21]{Lee} and \cite[Theorem 2.5]{Adachi}. By general position, this disk meets itself at most in isolated points. Sliding $h_{i}$ along this disk we may, by handle passes only, make $h_{i}$ into a handle connecting $B_{1}$ to $B_{i}$ without passing through any bases. Now all handles in $H$ may be moved by handle passes so they only pass through base $B_{1}$, after which all other bases may be deleted after appropriate handle slides. But a ribbon presentation with a single base is always stably equivalent to the unknotted ribbon presentation. \epr

For a ribbon $n$-knot $K$ with $n\geq 2$ and ribbon presentation $(B, H)$, there is a natural presentation for the knot quandle $Q(S^{n+2}, K)$. This presentation is generated by the bases in $B$ and has one relation for each handle $h_{j}\in H$. The relation $r_{j}$ has the following form: let $b_s$ and $b_e$ be the starting and ending bases for $h_j$, and let $h_j$ follow the word $w=b_1 b_2... b_n$ (here $w$ is a word in the bases in $B$, where each base can appear positively or negatively, depending on whether the handle passes through it in a direction that agrees with the normal vector or opposite to it) as it moves from $b_s$ to $b_e$. Then the relation $r_j$ is given by the equation $b_e=b_s^{\overline{b_1}\overline{b_2}...\overline{b_n}}$. Geometrically, the generator corresponding to a base $b\in B$ is given by an equivalence class of paths containing a path that goes from the basepoint to the base $b$ without passing through any other bases.

\begin{theorem}
A ribbon $n$-knot, $n\geq 2$, with ribbon presentation $(B, H)$, has a knot quandle whose presentation is as described in the previous paragraph.
\end{theorem}
\bpr
This is proved by repeated application of the Van Kampen theorem. Alternately, it is not difficult to write down the double point curves in a projection of a ribbon knot where each base and handle is embedded under the projection. Then the double point curves correspond to handles passing through bases, and the relations may be worked out using the computations as shown in e.g. \cite{CKS} (which are themselves based on the Van Kampen theorem).
\epr

\begin{remark}
Let $Q(S^{n+2}, K)$ be the quandle of a ribbon knot $K$ with presentation $(B, H)$. To specify an element of $Q(S^{n+2}, K)$, it suffices to specify a path $\gamma$ by assuming that $\gamma$ starts at the basepoint for the quandle, then specifying the ordered signed intersections of the path with the bases, and finally specifying on which base the path terminates.
\end{remark}

This follows from general position arguments.

\begin{remark}
We will denote the generator of a quandle corresponding to a base $b$ by using the same symbol; whether we are referring to the base or the quandle element will be clear from context.
\end{remark}

Throughout the following, we will make extensive use of the following fact, which forms part of the geometric definition of the quandle of a knot complement. Let $K$ be a knot in $S^{n+2}$, and let $x$ be the basepoint for the quandle $Q(S^{n+2}, K)$. Let $N(K)$ be the closure of a tubular neighborhood of $K$. Recall that elements of the quandle are equivalence classes of paths $\gamma$ which start on $x$ and terminate on $\partial N(K)$. Two paths $\gamma, \gamma'$ are equivalent as quandle elements iff they are homotopic through paths starting on $x$ and terminating on $\partial N(K)$. Therefore, there exists a disk $D_{\gamma, \gamma'}$ whose boundary consists of the union of $\gamma \cup \gamma'$ with some path on $\partial N(K)$. By various approximation theorems, \cite[Theorem 10.16, Theorem 10.21]{Lee} and \cite[Theorem 2.5]{Adachi}, we may assume this disk is smooth and immersed, and therefore by general position we may assume that $D_{\gamma, \gamma'}$ is immersed and embedded except at a finite number of points, provided that $n\geq 2$. In fact for $n\geq 3$, general position arguments show that the disk can be taken to be embedded. Thus, the homotopy may be taken to be a homotopy through embedded paths, that is, a smooth isotopy of the path which traces out the disk $D_{\gamma, \gamma'}$.

We will repeatedly make use of the existence of such isotopies. All our applications are essentially uses of the following construction. Let $(B, H)$ be a ribbon presentation for $K$. Let $h, h'$ be handles such that $(B, H\cup \{h\})$ and $(B, H\cup \{h'\})$ are ribbon presentations. Recall that a handle is an embedding of $D^{1}\times D^{n}$ into $S^{n+2}$ such that $\{0,1\}\times D^{n}$ is embedded in elements of $B$. We will refer to $h|_{[0,1]\times\{0\}}$ as the \emph{core} of $h$. Observe that once the core of $h$ is specified, the handle itself is specified up to isotopy. Now suppose that $h, h'$ have the property that $h|_{[0,1/2]\times\{0\}}=h'|_{[0,1/2]\times\{0\}}$. Let $\gamma=h|_{[1/2,1]\times\{0\}}, \gamma'=h'|_{[1/2,1]\times\{0\}}$. Now $\gamma, \gamma'$ determine quandle elements for $Q(S^{n+2}, K)$.

\begin{lemma}
Let $K$, $(B, H)$, $h, h', \gamma, \gamma'$ be as described in the previous paragraph, and suppose that $\gamma\cong\gamma'$ as elements of $Q(S^{n+2}, K)$. Then there is a stable equivalence between $(B, H\cup \{h\})$ and $(B, H\cup \{h'\})$ which involves only performing handle slides and passes isotoping $h$ to agree with $h'$.\label{qhs}
\end{lemma}
\bpr Let $D_{\gamma, \gamma'}$ be the homotopy disk giving the equivalence of $\gamma, \gamma'$ as elements of $Q(S^{n+2}, K)$. This is a disk in the complement of $K$. This disk may be perturbed so that its interior does not meets $h$ or $h'$, by general position arguments. Now $D_{\gamma, \gamma'}$ intersects itself only in isolated points. Therefore, we may slide $h|_{[1/2, 1]\times\{0\}}$ along this disk until it agrees with $h'|_{[1/2, 1]\times\{0\}}$. During this process we will be sliding the end of the core of $h$ around on $\partial N(K)$ so we may perform some handle slides, as well as passes. \epr

\begin{remark}
This argument cannot be extended to an argument for an analogous result for handles $h$ and $h'$ that are already part of the ribbon presentation for $K$. This is because in that case, the disk $D_{\gamma, \gamma'}$ might pass through the interior of the handle $h$ or $h'$ themselves. This cannot happen in the above case because $\gamma, \gamma'$ are homotopic using a homotopy in the complement of $K$, and, therefore, the homotopy cannot move the path through the interior of $h$ or $h'$ since they are not handles for the presentation of $K$.
\end{remark}

\begin{lemma}
Let $K$ be a ribbon $n$-knot, $n\geq 2$, with ribbon presentation $(B, H)$. Suppose that the quandle $Q(S^{n+2}, K)$ is generated by $B-\{ b\}$. Then, by adding a trivial handle to $(B, H)$, we can merge $b$ and some $b'\in B-\{ b\}$ into a single base through stable equivalence.\label{reduce}
\end{lemma}
\bpr Let $h$ be the new trivial handle, with both ends on $b$. We may isotope $h$ so that the midpoint of the core of $h$ passes through the basepoint $x$ for the quandle while it remains a trivial handle. Now the path that goes from the basepoint to $b$ is equivalent in the quandle to a path $\gamma$ from the basepoint to some other base $b'$ which does not pass through the base $b$, since $B-\{ b\}$ generates the quandle. There is a quandle homotopy disk connecting these two paths which does not intersect $h$ (since all the relations in the quandle are geometrically given by the handles in $H$). We may slide the portion of $h$ between its midpoint and terminal point along this disk so that it joins $b$ and $b'$ and follows the path $\gamma$, using the method described in Lemma \ref{qhs}. Since $\gamma$ does not pass through $b$, we can now perform handle passes through $h$ so that any handles passing through $b$ are moved to instead pass through $b'$. Once $b$ no longer contains any ribbon intersections, we may collapse one of the attaching handles to join it to an adjacent base. \epr


\begin{theorem}
Let $K, K'$ be ribbon $n$-knots, $n\geq 2$, with ribbon presentations $(B, H), (B, H')$, respectively, of the same genus, such that the identity map on $B$ induces a quandle isomorphism between $Q(S^{n+2}, K)$ and $Q(S^{n+2}, K')$. Then they are $|H|$-weakly stably equivalent. \label{isoeq}
\end{theorem}
\bpr Let $h'\in H'$ connect $b_0, b_1$. Attach a trivial handle $h$ to $H$ with both ends on $b_0$. Let $h$ be isotoped so that the midpoint of its core passes through the basepoint $x$ for the quandle. The path from $x$ to $b_0$ following the core of $h$ and the path from $x$ to $b_1$ following $h'$ are homotopic in the quandle $Q(S^{n+2}, K)$; thus we may slide $h$'s latter half along the homotopy disk, as described in Lemma \ref{qhs}. After this is done, $h$ will pass through exactly the same bases as $h'$. Repeat this for each handle in $H'$. Then perform the same operation using the trivial handles we have added to $H'$, modifying them so that all the handles originally in $H$ are now in $H'$ as well. We now have all the handles in $H$ in $H'$ as well, and all the handles in $H'$ have been added to $H$. Thus we have modified the two ribbon presentations so they have identical sets of bases and handles. \epr

\begin{theorem}
Let $K, K'$ be ribbon $n$-knots, $n\geq 2$, with ribbon presentations of the same genus $(B, H), (B', H')$, respectively. Let $f:B\rightarrow Q(S^{n+2}, K')$ be a function which extends to a quandle isomorphism $Q(S^{n+2}, K)\rightarrow Q(S^{n+2}, K')$. Then $K, K'$ are $(|B'|+|H|)$-weakly stably equivalent with bases given by $B$ and $f(B)$ respectively. \label{injeq}
\end{theorem}
\bpr For each $b_i\in B$, we obtain a quandle element $b'_i=f(b_i)\in Q(S^{n+2}, K')$. The $b_i'$ are paths from the basepoint to some base $B'_i \in B'$. Create a trivial base $B''_i$ connected by a handle to $B'_i$ by a handle that does not pass through any bases. Then we may slide $b_i'$ so that it terminates on $B''_i$. Clearly $b_i'$ does not pass through $B''_i$. Therefore, we can pull $B''_i$ back along $b_i'$, contracting the quandle element, until we have so isotoped $B''_i$ such that $b_i'$ is now a path from the basepoint to $B''_i$ which does not pass through any bases.

Now the collection of $B''_i$s generates the quandle $Q(S^{n+2}, K')$. By Lemma \ref{reduce}, $(B', H')$ is $|B'|$-weakly stably equivalent to a knot whose bases consist of exactly the $B''_i$s. Add an equal number of handles to $H$. Then by Theorem \ref{isoeq} the result follows. \epr

\begin{theorem}
Let $K, K'$ be ribbon $n$-knots (of the same genus) with ribbon presentations $(B, H), (B', H')$, respectively. Suppose the knot quandles $Q(S^{n+2}, K)$ and $Q(S^{n+2}, K')$ are isomorphic. Then $K$ and $K'$ are $(|B'|+|H|)$-weakly stably equivalent.
\end{theorem}
\bpr Let $f:Q(S^{n+2}, K)\rightarrow Q(S^{n+2}, K')$ be the isomorphism. Then $f$ defines a map from $B$ to $Q(S^{n+2}, K')$ defined by sending $b\mapsto f(b)$ (where we send the base to the quandle element which its standard generator is sent to by $f$. Now apply Theorem \ref{injeq}. \epr

\begin{corollary}
Let $K, K'$ be as in the previous theorem, but with possibly different genus. Then they are weakly stably equivalent.
\end{corollary}

\begin{corollary}
Let $(B, H), (B', H')$ be two ribbon presentations for a single $n$-knot $K$, $n\geq 2$. Then they are $(|B'|+|H|)$-weakly stably equivalent.
\end{corollary}

This is related to the question posed in \cite{Naka}, where it is asked whether two ribbon presentations for the same $n$-knot with $n\geq 2$ are necessarily stably equivalent. In \cite{NakaNaka} an infinite family of examples are constructed showing that two ribbon presentations for the same $1$-knot do not need to be stably equivalent. In fact examples are constructed of $1$-knots with $n$ stably inequivalent ribbon presentations for any $n\in\mathbb{N}$.

\begin{remark}
These results all hold in their natural generalizations to ribbon links with multiple components.
\end{remark}

%

\end{document}